\documentclass[a4paper,10pt]{article}
\usepackage{amsfonts}
\usepackage{amssymb}
\usepackage{amsmath}
\usepackage{cite}
\newtheorem{theo}{Theorem}[section]
\newtheorem{lema}[theo]{Lemma}
\newtheorem{rema}[theo]{Remark}
\newtheorem{prop}[theo]{Proposition}

\newtheorem{coro}[theo]{Corolary}
\newtheorem{defi}[theo]{Definition}	
\newcommand{\B}{\mathbb{B}}

\newcommand{\R}{\mathbb{R}}
\newcommand{\D}{\mathbb{D}}
\newcommand{\Td}{\mathbb{T}^d}
\newcommand{\T}{\mathbb{T}^1}

\newcommand{\C}{\mathbb{C}}
\newcommand{\N}{\mathbb{N}}
\newcommand{\eps}{\varepsilon}

\newcommand{\tal}{\theta+\alpha}

\newcommand{\Z}{\mathbb{Z}}

\newcommand{\Pru}{\mathbf{Proof. }}

\begin{document}
\title{Local dynamics for fibered holomorphic transformations}
\author{Mario PONCE\\ Universit\'e Paris SUD-XI}
\maketitle
\begin{abstract}
Fibered holomorphic dynamics are skew-product transformations $F(\theta,z)=(\tal, f_{\theta}(z))$ over an irrational rotation, whose fibers are holomorphic functions. In this paper we study such a  dynamics on a neighborhood of an invariant curve.  We obtain some results  analogous to the results  in the non fibered case. 
\end{abstract}
\section{Introduction}
Let $\alpha$ be a rationally independent vector in $\Td$, which will be  fixed for the rest of this work. Let $U\subset \C$ be an open simply connected neighborhood of the origin. We consider fibered continuous injective transformations  
\begin{eqnarray*} 
F:\Td\times U &\longrightarrow& \Td\times \C\\
 (\theta,z)&\longmapsto& \big( \tal,f_{\theta}(z)\big).
 \end{eqnarray*}
We will assume that the functions $f_{\theta}:U\to \C$ are holomorphic for all $\theta\in \Td$; we call $F$ a fibered holomorphic dynamics, and we denote it by  \emph{fhd}.

The \emph{fhd}'s are a special class of the well-known \emph{skew-product} transformations. A closely related class of skew-product transformations over irrational rotations are fibered circle homeomorphisms, where the fiber is a circle and the $f_{\theta}$ are circle homeomorphisms . In \cite{HERM83} M.Herman establishes    the  basis for the study of the fibered circle homeomorphisms, defining the fibered rotation number. Recently this study has been relaunched mainly by the works of G.Keller, J.Stark and T.J\"ager (see \cite{JAKE06},\cite{JAST06})). They have established in particular a Poincar\'e-like classification of these transformations by means of the fibered rotation number and the presence of invariant graphs. In his doctoral thesis, O.Sester \cite{SEST97} (see also \cite{SEST99}) has studied hyperbolic fibered polynomials,  successfully generalizing the classical notions of Julia set, Green's function and the principal cardioid of the Mandelbrot set in the parameter space. Other important contributions to this subject are the works of M.Jonsson \cite{JONS99}, \cite{JONS00}. 

 The notion of periodic or fixed points for $F$ has no sense, since the irrational rotation in the base is minimal. We are interested in the local dynamics of $F$ in a neighborhood of an invariant curve, that is,  a continuous curve $u:\Td\to U$ such that
\begin{equation}\label{introequinv}
 F\big(\theta,u(\theta)\big)=\big(\theta+\alpha,u(\theta+\alpha)\big)
 \end{equation}
 for all $\theta\in \Td$, or equivalently $f_{\theta}\big(u(\theta)\big)=u(\tal)$. These objects play the role of a center around which  the dynamics of $F$ is organized, generalizing thus the role of a fixed point for the local dynamics of an holomorphic germ $g:(\C,0)\to(\C,0)$, the so-called \emph{non fibered case}. The rest of this work is devoted to precise and prove this last assertion. A complete survey on the non fibered case can be found in    \cite{CAGA93}.

\section{Definitions}
Let $F$ be a \emph{fhd} and $u$ an invariant curve for $F$. We define the infinitesimal characteristics of the curve $u$ by
\begin{defi} The \emph{fibered multiplicator} $\kappa(u)$ of the curve is the real number
\begin{equation*}
\kappa(u)=exp\Big(\int_{\Td}\log \big|\partial_zf_{\theta}\big(u(\theta)\big)\big|d\theta\Big).
\end{equation*} 
We recall that $F$ is injective and so the differential $\partial_zf_{\theta}$ is always non zero. We say $u$ is an attracting curve if $\kappa(u)<1$,  $u$ is a repelling curve if $\kappa(u)>1$ and an indifferent curve if $\kappa(u)=1$. 
\end{defi}
%\begin{defi}
%The \emph{degree} $deg(u)$ of the curve is the integer number corresponding to the topological degree of the application $\theta\mapsto \partial_zf_{\theta}\big(u(\theta)\big)$. Note that the application $\theta\mapsto exp\big(-2\pi i deg(u)\theta\big)\partial_zf_{\theta}\big(u(\theta)\big)$ has zero topological degree and then we can define its logarithm. 
%\end{defi}
\begin{defi}
%Let $\kappa=\kappa(u)$ and $n=\deg(u)$. The number $\varrho_{tr}(u)\in \T$ is defined so that
%\begin{equation*}
%\int_{\T}\log \Big(exp\big(-2\pi i n\theta\big)\partial_zf_{\theta}\big(u(\theta)\big)\Big)d\theta=\log \kappa+2\pi i \varrho_{tr}(u).
%\end{equation*}
Suppose $\kappa(u)=1$ and the  application 
\begin{equation*}
\theta\longmapsto \partial_{z}f_{\theta}\big(u(\theta)\big)
\end{equation*}
 is homotopic  in $\C\setminus \{0\}$ to a constant. We refer to this situation as the \emph{indifferent zero degree case}. We define a number which represents  the average rotation speed of the dynamics around the invariant curve, 
\begin{equation*}
\varrho_{tr}(u)=\frac{1}{2\pi i }\int_{\Td}\log \partial_zf_{\theta}\big(u(\theta)\big)d\theta .
\end{equation*}
We call this number  the \emph{fibered rotation number}.  Note the $\log$ above is well defined $\mod 2\pi i$ and the number $\varrho_{tr}(u)$ is well defined $\mod 1$. 
\end{defi}

In general, we will  consider the invariant curve as being the zero section curve $\{z\equiv 0\}_{\Td}$. We obtain this situation by conjugating $F$ with the change of coordinates given by $H(\theta)=\big(\theta,z+u(\theta)\big)$. In this way the resulting transformation $\tilde{F}=H^{-1}\circ F\circ H$ has the zero section as an invariant curve with the same infinitesimal characteristics as the original curve $u$ for $F$. More generally,  we will need  to consider continuous change of coordinates $\tilde{H}$  defined on a tubular neighborhood of the invariant curve (the zero section), say
\begin{equation*}
(\theta,z)\stackrel{\tilde{H}}{\longmapsto}\big(\theta,h_{\theta}(z)\big).
\end{equation*}
The functions $h_{\theta}$ will be  biholomorphic transformations between two topological discs fixing the origin, that is $h_{\theta}(0)=0$. When conjugating our transformation $F$ we will get a new transformation $\tilde{F}=H^{-1}\circ F\circ H$ having the zero section $\tilde{u}=\{z\equiv 0\}_{\Td}$ as an invariant curve and $\kappa(\tilde{u})=\kappa(u)$. That is, the fibered multiplicator is invariant by conjugacy. In the indifferent zero degree case we will admit only those changes of coordinates which have themselves  zero degree, that is, the application $\theta\mapsto h_{\theta}(0)$ is homotopic in $\C\setminus \{0\}$ to a constant.  In that case  we get also $\varrho_{tr}(\tilde{u})=\varrho_{tr}(u)$. In the indifferent zero degree case the fibered rotation number is also invariant by zero degree conjugacy. 

We say $F$ is linearizable if we can conjugate it on a tubular neighborhood of the invariant curve to a transformation $\Lambda_A$ of the form
\begin{equation*}
\Lambda_A(\theta,z)=\big(\tal,A(\theta)z\big),
\end{equation*} 
for a continuous complex function $A:\Td\to \C$. When the absolute value $|A(\theta)|=\kappa(u)$ for all $\theta\in \Td$ we say $F$ is \emph{modulus linearizable}. When $|A(\theta)|<1$ for every $\theta\in \Td$ ($|A(\theta)|>1$ for every $\theta\in \Td$) we say that $F$ is \emph{weakly linearizable} to an attractive (repelling) linear \emph{fhd}. In the indifferent zero degree case we say $F$ is \emph{strongly linearizable} if $A(\theta)=e^{2\pi i \varrho_{tr}(u)}$.

We  say an invariant curve $u$ is \emph{stable} if there exists an open invariant and bounded tubular neighborhood $\mathcal {U}$ , $F(\mathcal{U})\subset \mathcal{U}$.  We always suppose that the  fibers $\mathcal{U}_{\theta}$ are topological discs and that $\mathcal{U}$ contains the curve. We call those neighborhoods an \emph{open invariant tube}.  
\section{Statement of the results}
In this paper we prove the following results concerning the local dynamics of a \emph{fhd} $F$ around an invariant curve  $u$:
\begin{prop}\label{prop1}
If $\kappa(u)<1$ then there exist an open tubular neighborhood  of the curve which is attracted to the curve by positive iteration. Moreover, $F$ is weakly linearizable.
\end{prop}
 When $\kappa(u)>1$ we obtain the analogous result for the repelling case by considering the inverse $F^{-1}$. When the curve is indifferent ($\kappa(u)=1$) we generalize to fibered dynamics a very well-known fact in the non fibered case
 \begin{center}
 \emph{Lyapounov stability is equivalent to linearization.}
 \end{center}

\begin{prop}\label{prop2}
If the curve is indifferent then stability of the curve is equivalent to modulus linearization.
\end{prop}
As we will see, linearizability alone does not implies stability.
 Then we  look at non stable situation in the indifferent case. We show that there are still some nearby points with a complete  orbit (past and future) which stays near the invariant curve.   Indeed, we prove a fibered version of the Continua's Theorem by R.P\'erez-Marco \cite{PERE97}
 \begin{theo}[Fibered P\'erez-Marco's Continua]\label{teo1}
 Let $F$ be a \emph{fhd} with an indifferent invariant curve $u$. Let $\mathcal{U}$ be an open neighborhood of the curve whose fibers are Jordan domains (the interior of a Jordan curve). We also assume  that $\overline{\mathcal{U}_{\theta}}$ depends continuously on $\theta$ and that $F$, $F^{-1}$ define \emph{fhd}'s which are injective in some neighborhood of $\overline{\mathcal{U}}$. Then there exists a compact connected set $K\subset \overline{\mathcal{U}}$ such that
 \begin{itemize}
 \item[$i)$ ] Fibers $K_{\theta}$ are connected full compact sets for every $\theta\in \Td$.
 \item[$ii)$ ] $graph(u)\subset K$.
 \item[$iii)$ ] $F(K)=F^{-1}(K)=K$, that is, $K$ is completely invariant by $F$.
 \item[$iv)$ ] $K\cap\partial\mathcal{U}\neq \emptyset$.
 \end{itemize}
  \end{theo}

  We also consider \emph{fhd}'s which are analytic with respect to the $\theta$ variable. Let $\delta$ be a positive real number. An analytic \emph{fhd} is a transformation 
$
  (\theta,z)\longmapsto F(\theta,z)
$
  that is defined and holomorphic as a $d+1$ variables function on the product $B_{\delta}^d\times U$, where  $B_{\delta}=\{\theta\in \C/\Z \ | \ |Im (\theta)|<\delta\}$. For this class of regularity and under the indifferent zero degree hypothesis  one can show  a fibered version of the Siegel's linearisation theorem \cite{SIEG42}. The Siegel's theorem for holomorphic germs  state that under some arithmetical condition on the rotation number,  the dynamics  is linearizable without additional hypothesis on the stability of the fixed point. More precisely one has
  \begin{theo}[Siegel's theorem for \emph{fhd}]
  Let $\delta>0$, $F$ be a \emph{fhd} analytic in the product $B_{\delta}^d\times \D$. Let $u:\B_{\delta}^d\to \C$ be an indifferent invariant analytic curve with zero degree. If the fibered rotation number $\varrho_{tr}(u)=\beta$ is such that the pair $(\alpha,\beta)$ verifies the arithmetical condition $\C\D_{\geq 1}(c,\tau)$ for some $c>0, \tau\geq 0$, then $F$ is strongly linearizable on a tubular neighborhood of the curve.
  \end{theo}
  The arithmetical condition $\C\D_{\geq 1}(c,\tau)$ is defined by
  \begin{equation*}
  \C\D_{\geq 1}(c,\tau)=\Big\{(\alpha,\beta)\in \Td\times \T \ \big | \ \|n\cdot\alpha-j\beta\|>\frac{c}{(|n|+j)^{2+\tau}}\ \forall\ j\geq 0, n\in \Z^d\Big\},
  \end{equation*}
  where for an integer vector $n=(n_1,\dots,n_d)\in \Z^d$ we put $|n|=\max_{1\leq i\leq d}|n_i|$.
    The proof  of this theorem is a straightforward adaptation of the classical proof due to Moser  for the Siegel's theorem (see \cite{KATO95}), and can be found in \cite{PONC07t}.
 \section{Proofs} 
 We recall that the invariant curve is supposed to be the zero section curve $u=\{z\equiv 0\}_{\T}$.
 We use the standard notation in the skew-product dynamics
 \begin{equation*}
 f_{\theta}^n(z):=\Pi_2F^n(\theta,z)=f_{\theta+(n-1)\alpha}\Big(f_{\theta+(n-2)\alpha}\big(\dots(f_{\theta}(z))\dots\big)\Big)
 \end{equation*}
 for the holomorphic second coordinate. We put $\rho_1(\theta)=\partial_zf_{\theta}(0)$, and in this way the transformation $F$ has the form
 \begin{equation}\label{notacion}
 F(\theta,z)=\Big(\tal, \rho_1(\theta)z+z^2\rho(\theta,z)\Big).
 \end{equation}
 The function $\rho$ is a continuos function, holomorphic in each fiber and $\rho_1:\Td\to \C\setminus\{0\}$ is a continuous function.  
 \subsection{Proof of proposition \ref{prop1}}
 Previous to the main proof, we want  to point at the following simple complex analysis result
  \begin{lema}\label{lemainversion}
 For any constants $A\in (0,1), B>0$ there exists constants $R=R(A,B),r=r(A,B)\in (0,1)$  such  that if $f:\D\to \C$ is an holomorphic function verifying
 \begin{itemize}
 \item[$i) $] $f(0)=0$,
 \item[$ii) $] $\|f\|_{\D}:=\sum_{\D} |f(z)|\leq B$,
 \item[$iii) $] $A<\big|\partial_zf(0)\big|<A^{-1}$,
 \end{itemize}
 then  there exists an open set $U_f\supset D(0,r)$ such that $f:U_f\to D(0,R)$ is a biholomorphism $\quad_{\blacksquare}$
 \end{lema}
As a corollary we obtain
 \begin{coro} Let $F:\Td\times \D\to \Td\times \C$ be  a \emph{fhd}, $\alpha\in \Td$ its rotation number on the base $\Td$. Let $u:\Td\to \D$ be an invariant curve whose multiplicator is  $\kappa=\kappa(u)$. Then there exists a positive radius  $R\in (0,1)$ such that the inverse transformation  $F^{-1}$ is a  \emph{fhd} which is well defined in the tube $\mathcal{U}_R$, of radius  $R$ and centered on the invariant curve. The rotation number of $F^{-1}$ on the base is  $-\alpha$. The curve $u$ is invariant by $F^{-1}$ and its multiplicator as an invariant curve for $F^{-1}$ is  $\kappa^{-1}\quad_{\blacksquare}$
 \end{coro}

We now begin the proof of proposition \ref{prop1}. 
 We assume $\int_{\Td}\log \big|\rho_1(\theta)\big|d\theta<0$ and we put $\kappa=\kappa(u)\in (0,1)$. We introduce a first conjugacy in order that the absolute value  $|\rho_1(\theta)|$ becomes closer to his mean $\kappa$. Let $l:\Td\to \R$ be a trigonometric polynomial, $l(\theta)=\sum_{0\leq |n|<N}e^{2\pi i n(\cdot\theta)}\hat{l}(n)$, verifying
 \begin{enumerate}
 \item[$i)$ ] $\hat{l}(0)=\int_{\Td}l(\theta)d\theta=\log \kappa$.
 \item[$ii)$ ] $\Big|\log \big|\rho_1(\theta)\big|-l(\theta)\Big|<\log \frac{1+\kappa^{1/2}}{2}$ for all $\theta\in \Td$.
 \end{enumerate} 
Let $\tilde{u}_1:\Td\to \R$ be the zero-mean solution to the cohomological equation 
 \begin{equation*}
 \tilde{u}_1(\theta)-\tilde{u}_1(\tal)=\log \kappa -l(\theta).
 \end{equation*}
The Fourier series method gives that $\tilde{u}_1(\theta)=\sum_{0<|n|<N}\frac{\hat{l}(n)}{e^{2\pi i (n\cdot\alpha)}-1}e^{2\pi i n\theta}$. We put $u_1(\theta)=e^{\tilde{u}_1(\theta)}$. We conjugate $F$ by the fibered re-scaling $H(\theta,z)=\big(\theta,u_1(\theta)z\big)$ getting 
\begin{eqnarray*}
  \big|\partial_z\big(H^{-1}\circ F\circ H\big)_{\theta}(0)\big|&=&\frac{u_1(\theta)}{u_1(\tal)}|\rho_1(\theta)|\\
  &=&\kappa \frac{|\rho_1(\theta)|}{e^{l(\theta)}}<\frac{\kappa+\kappa^{1/2}}{2}<1.
  \end{eqnarray*}
We thus may suppose that $|\rho_1(\theta)|<\frac{\kappa+\kappa^{1/2}}{2}<1$ for all $\theta\in \Td$, and the holomorphic functions $f_{\theta}$ are defined in a disc $D_R$ for some $R>0$. The following proof is an adaptation of the classical proof of the non fibered corresponding result (see \cite{CAGA93}).\\
We define a continuous  function $\phi:\Td\times D_R\to \C$ holomorphic in each fiber, by the equality
  \begin{equation*}
  f_{\theta}(z)=\rho_1(\theta)z\big\{1+\phi(\theta,z)z\big\}.
  \end{equation*}
  Let $\{g^n:\Td\times D_R\to \C\}_{n\geq 1}$ be the sequence of continuous functions  defined by
  \begin{equation*}
  g_{\theta}^n(z)=\Big(\prod_{i=0}^{n-1}\rho_1(\theta+i\alpha)\Big)^{-1}f_{\theta}^n(z).
  \end{equation*}
  One has
  \begin{eqnarray*}
  f_{\theta}^{n+1}(z)&=&\rho_1(\theta+n\alpha)f_{\theta}^n(z)
\big\{1+\phi\big(\theta+n\alpha,f_{\theta}^n(z)\big)f_{\theta}^n(z)\big\}\\
g_{\theta}^{n+1}(z)&=&g_{\theta}^n(z)\big\{1+\phi\big(\theta+n\alpha,f_{\theta}^n(z)\big)f_{\theta}^n(z)\big\}\\
&=&z\prod_{j=0}^n\big\{1+\phi\big(\theta+j\alpha,f_{\theta}^j(z)\big)f_{\theta}^j(z)\big\}.
  \end{eqnarray*}
  We define another  continuous function $\psi:\Td\times D_R\to\C$ by
  \begin{equation*}
  \log \big(1+\phi(\theta,z)z\big)=z\psi(\theta,z).
  \end{equation*}
The function $\psi$ is bounded. We also observe that the $g^n$ and $\psi$ are holomorphic in each fiber. We want  to show that the sequence $\{g^n\}_{n\geq 0}$ is uniformly convergent on a product $\Td\times D_r$, $0<r\leq R$ to a (continuous, holomorphic in each fiber) function $g:\Td\times D_R\to \C$. To do that is sufficient  to show that the following series are an absolutely and uniformly convergent  series
  \begin{equation}\label{stable.seriepsi}
  \sum_{j=0}^{\infty} \psi\big(\theta+j\alpha,f_{\theta}^j(z)\big)f_{\theta}^j(z).
  \end{equation}
From  the bound $\big|\rho_1(\theta)\big|<\frac{\kappa+\kappa^{1/2}}{2}<1$ we get a constant $c\in (0,1)$ and a positive radius $r$ such that $\big|f_{\theta}^n(z)\big|<c^n|z|$ for $|z|<r$, $\theta\in \Td$ and for all $n$ large enough. This implies the desired convergence. Thus, the limit function $g:\Td\times D_r\to \C$ is continuous, holomorphic in each fiber and verifies $g_{\theta}(0)=0$, $\partial_zg_{\theta}(0)=1$ for all $\theta\in \Td$. Lemma \ref{lemainversion} applies and give us a tubular neighborhood of the zero section where $g$ is invertible. There, the inverse $g^{-1}$ is a continuous function; defining
  \begin{eqnarray*}
  \Phi_n(\theta,z)&=&\big(\theta,g^n_{\theta}(z)\big)\\
  \Phi(\theta,z)&=&\big(\theta,g_{\theta}(z)\big),
  \end{eqnarray*}
we have 
  \begin{eqnarray*}
  \Phi_n\circ F&=&\Big(\tal,\Big(\prod_{i=0}^{n-1}\rho_1(\theta+\alpha+i\alpha)\Big)^{-1}f^n_{\tal}\big(f_{\theta}(z)\big)\Big)\\
  &=&\Big(\tal,\rho_1(\theta)\Big(\prod_{i=0}^{n}\rho_1(\theta+i\alpha)\Big)^{-1}f_{\theta}^{n+1}(z)\Big)\\
  &=&\Lambda_{\rho_1}\circ \Phi_{n+1}.
  \end{eqnarray*}
 Thus in the limit 
  \begin{equation*}
   F=\Phi^{-1}\circ \Lambda_{\rho_1}\circ \Phi\quad_{\blacksquare}
  \end{equation*}

  The corresponding result in the repelling case ($\kappa(u)>1$) is an easy consequence of the attractive one, by considering the curve as an invariant attractive curve for the inverse $F^{-1}$, which is a well defined \emph{fhd} on a tubular neighborhood of the curve.
  \subsection{Proof of  proposition \ref{prop2}}
  If $F$ is modulus linearizable then an open invariant tube is found by the image of $\Td\times D_r$ by the conjugacy function,  for a small $r>0$. We must   prove then the direct implication. We will call \emph{domain} to any simply connected open set in the complex plane $\C$, containing the origin and different of  $\C$.
  We need to recall some definitions about domains  and their uniformizations. We refer the reader  to the book of C.Pommerenke \cite{POMM75} for further references about the facts below.
  \paragraph{Conformal Radius. } Let $\Omega$ be a domain and let $h:\D\to \Omega$ be the uniformization (biholomorphic function) verifying $h(0)=0$ and $h'(0)>0$. We call the positive real number $h'(0)$ the \emph{conformal radius} of $\Omega$ and we denote it by $R(\Omega)$. The following monotonicity property holds: If $\Omega'\subset \Omega$ are domains, then one has $R(\Omega')\leq R(\Omega)$ and equality occurs only in the case $\Omega'=\Omega$.  We associate to a domain $\Omega$ a sequence $\{\Omega^n\}_{n\geq 2}$ of domains defined by $\Omega^n=h\big(D_{1-\frac{1}{n}}\big)$. This sequence is a monotone exhausting sequence  and the conformal radius verify $R(\Omega^n)\nearrow R(\Omega)$.
  \paragraph{Caratheodory's Kernel. } Let $\{O_n\}_{n\geq 0}$ be a sequence of domains. The Kernel $\nabla(\{O_n\})$ of this sequence consists of the union of the origin and those points $z\in \C$ which satisfy the following property: There exists a domain $\Omega_z$ containing $z$, such that $\Omega_z\subset O_n$ for $n$ large enough. The kernel is either reduced to the origin, a domain, or $\C$.
  \paragraph{Kernel convergence. } We say the sequence $\{O_n\}_{n\geq 0}$ is kernel convergent to a domain $O$ if every subsequence of $\{O_n\}_{n\geq 0}$ has the domain $O$ as Kernel.
  \begin{theo}[Caratheodory]\label{caratheodory} Let $\{f_n:\D\to \C\}_{n\geq 0}$ be a sequence of univalent functions, with $f_n(0)=0$ and $f_n'(0)>0$ for all $n\geq 0$. Let  $O_n=f_n(\D)$ and $O$ be domains. Then the sequence $\{O_n\}$ is kernel convergent to $O$ iff the sequence $\{f_n\}$ is convergent  uniformly on compacts of $\D$. Moreover, in this case the limit function $f$ is a conformal representation of $O$. Finally, if $K\subset O$ is a compact set, then $K$ is contained in $O_n$ for large $n$ and $f_n^{-1}\big|_K$ converges uniformly to $f^{-1}\big|_K$.
  \end{theo}
  
  For $n\geq 0$ and $\theta \in \Td$ we define the Birkhoff sums associated to the invariant curve as being
   \begin{equation*}
  S_u^n(\theta)=\sum_{i=0}^{n-1}\log \big|\partial_zf_{\theta+i\alpha}\big(u(\theta+i\alpha)\big)\big|.
  \end{equation*}
In the indifferent case we have the uniform limit
\begin{equation*}
\log \kappa(u)=0=\lim_{n\to \infty}\frac{1}{n}S_u^n(\theta)
\end{equation*}
due to the unique ergodicity of the irrational rotation $\theta\mapsto \tal$. Let $\mathcal{U}$ be an invariant open tube containing the invariant curve (the zero section). Each fiber $\mathcal{U}_{\theta}$ is a domain. We associate to the each fiber $\mathcal{U}_{\theta}$ the uniformization $h_{\theta}:\D\to \mathcal{U}_{\theta}$ as before. 
\begin{lema} \label{lemasumas}There exist constants $B_+>0$,  $B_-<0$ and a point  $\tilde{\theta}\in \Td$ such that $B_-<S_u^n(\tilde{\theta})<B_+$ for all $n\geq 0$. 
\end{lema}
$\Pru $
Since the open set $\mathcal{U}$ is bounded, there exists a positive constant $W>0$ such that $W^{-1}<R\big(\mathcal{U}_{\theta}\big)<W$ for $\theta\in \T$. Let ${\theta}\in \T$ and $n\geq 0$. We define the functions $\tilde{f}_{\theta}^n:\D\to \D$ by
\begin{equation*}
\tilde{f}_{\theta}^n(z)=h^{-1}_{\theta+n\alpha}\Big(f_{\theta}^n\big(h_{\theta}(z)\big)\Big).
\end{equation*}  
The Schwartz Lemma yields $\big|\partial_z\tilde{f}_{\theta}^n(0)\big|\leq 1$, or in other words
\begin{equation*}
h'_{\theta+n\alpha}(0)^{-1}\Big(\prod_{i=0}^{n-1}|\rho_1(\theta+i\alpha)|\Big)h'_{\theta}(0)\leq 1.
\end{equation*} 
From this inequality the existence of $B_+$ is evident for any choice of $\tilde{\theta}$. We show now the existence of $\tilde{\theta}$ and $B_-$. We  proceed by contradiction. Let $l<0$. There would  exist a smallest natural number $n_l(\theta)$ verifying $S_u^{n_l(\theta)}(\theta)<l$.  There exists an uniform upper bound $N_l$ for $n_l(\theta)$. Indeed, otherwise,  there should exist a sequence $\{\theta_j\}_{j\in \N}\subset \Td$ such that $n_l(\theta_j)>j$ for all $j\geq 1$. Let $\hat{\theta}$ be an accumulation point of this sequence. Since the Birkhoff's sums are continuous functions in the $\theta$ variable, we know that in an open neighborhood of   that point  $\hat{\theta}$  one has $S_u^{n_l(\hat{\theta})}(\theta)<l$, a contradiction. 
\\
We put $l=-2B_+, N=N_l, n(\theta)=n_l(\theta)$. Then we have
\begin{equation*}
S_u^N(\theta)=S_u^{n(\theta)}(\theta)+S_u^{N-n(\theta)}(\theta+n(\theta)\alpha)<l+B_+=\frac{l}{2}<0
\end{equation*} 
but $\int_{\T}S_u^N(\theta)d\theta=0$, a contradiction $\quad_{\blacksquare}$
\\

Bounded Birkhoff's sums are an important tool when solving cohomological equations, as showed by the well-known 
 \begin{theo}[Gottschalk-Hedlund, see \cite{KATO95}]\label{gothed}
Assume $\mathcal{X}$ is a compact metric space, $f:\mathcal{X}\to \mathcal{X}$ is a continuous minimal transformation. Let $g:\mathcal{X}\to \C$ be a continuous function such that their  Birkhoff's sums verify
\begin{equation*}
\sup_{n\in \N}\Big|\sum_{i=0}^{n}g\circ f^i(x_0)\Big|<\infty
\end{equation*}
for some  $x_0\in \mathcal{X}$. Then there exists a continuous solution $h:\mathcal{X}\to \C$ for  the cohomological  equation
\begin{equation*}
h\circ f-h=g.
\end{equation*}
\end{theo}
\begin{rema}
As shown by Furstenberg,  there exists analytic minimal non-ergodic diffeomorphims of $\mathbb{T}^2$. The key of this result is the existence of a real irrational number $\omega$ and an analytic function $\phi:\T\to \R$ such that $\int_{\T}\phi d\theta =0$ and the cohomological equation
\begin{equation}\label{furstenberg}
\psi(\theta+\omega)-\psi(\theta)=\phi(\theta)
\end{equation} 
has measurable but not continuous solution (see \cite{KATO95}). Let $F(\theta,z)=\big(\theta+\omega, e^{\phi(\theta)}z\big)$. For this linear \emph{fdh} the invariant curve $\{z\equiv 0\}_{\T}$ is indifferent but not stable. Otherwise, lema \ref{lemasumas} and theorem \ref{gothed} provide a continuous solution to (\ref{furstenberg}).
\end{rema}

 Lemma \ref{lemasumas} allows us to apply the Gottschalk-Hedlund Theorem to the cohomological equation 
\begin{equation*}
\tilde{u}_1(\tal)-\tilde{u}_1(\theta)=\log \big|\rho_1(\theta)\big|.
\end{equation*}
We obtain a zero-mean continuous solution $\tilde{u}_1:\Td\to \R$. With $u_1(\theta)=exp\big(\tilde{u}_1(\theta)\big)$ and the change of variables $H(\theta,z)=\big(\theta,u_1(\theta)z \big)$ one has that the  conjugacy $\tilde{F}=H^{-1}\circ F\circ H$ verifies $\big|\partial_z\tilde{f}_{\theta}(0)\big|=1$ for all $\theta\in \Td$. The open tube $\tilde{\mathcal{U}}=H^{-1}\big(\mathcal{U}\big)$  is still an invariant open tube for $\tilde{F}$. Thus we may suppose that $|\rho_1(\theta)|=1$ for every $\theta \in \Td$. We put $R(\theta)=R\big(\mathcal{U}_{\theta}\big)$. 
    \begin{lema}\label{introRsci}
    The function $R:\Td\to \R^+$ is a lower semi continuous function.
    \end{lema}
    $\Pru$ Let $\eps>0$ and $\tilde{\theta}$ in $\Td$. Let $n^*\geq 2$ such that   $R(\Omega_{n^*}(\mathcal{U}_{\tilde{\theta}}))>R(\tilde{\theta})-\eps$. The set  $\mathcal{U}$ is open, therefore there exists an open neighborhood   $I_{\tilde{\theta}}$ around $\tilde{\theta}$ such that for  $s$ in such a neighborhood, the open set  $\mathcal{U}_s$ contains the compact set $\overline{\Omega_{n^*}(\mathcal{U}_{\tilde{\theta}})}$ (since this compact set is far away from the complementary set of  $\mathcal{U}$). One has then
    \begin{equation*}
    R(s)>R(\Omega_{n^*}(\mathcal{U}_{\tilde{\theta}}))>R({\tilde{\theta}})-\eps \quad_{\blacksquare}
    \end{equation*}
    \begin{lema}\label{introRconstante}
    The conformal radius  $R(\theta)$ is independent of  $\theta$ in $\Td$.
    \end{lema}
    $\Pru$  The function $R:\Td\to \R^+$ is l.s.c., thus there exists a point  $\tilde{\theta}$ in $\Td$ where $R$ reaches his minimum. From the inclusion
    \begin{equation*}
    f_{\tilde{\theta}-\alpha}\big(\mathcal{U}_{\tilde{\theta}-\alpha}\big)\subset \mathcal{U}_{\tilde{\theta}}
    \end{equation*}
we deduce that  
    \begin{equation*}
    R\Big(f_{\tilde{\theta}-\alpha}\big(\mathcal{U}_{\tilde{\theta}-\alpha}\big)\Big)\leq R(\tilde{\theta}).
\end{equation*}
The function $f_{\theta}$ preserves  the conformal radius since  $|\rho_1(\theta)|=1$ for every $\theta$ in $\Td$. Then  we have
 \begin{equation*}
 R\Big(f_{\tilde{\theta}-\alpha}\big(\mathcal{U}_{\tilde{\theta}-\alpha}\big)\Big)=R(\tilde{\theta}-\alpha)
 \end{equation*} 
 and the inequality 
 \begin{equation*}
 R(\tilde{\theta}-\alpha)=R(\tilde{\theta})
 \end{equation*}
from the minimality of the value $R(\tilde{\theta})$. The same argument tells us  that $R(\tilde{\theta}-n\alpha)=R(\tilde{\theta})$ for all $n\geq 0$. We have a l.s.c. function with a dense set of minimal points, therefore a constant function $\quad_{\blacksquare}$
 \\
 
 We define $R_F^{\mathcal{U}}$ as the conformal radius of some (any)  domain   $\mathcal{U}_{\theta}$.
 \begin{lema}\label{stable.232}
 Let $\{{\theta_{n}}\}_{n\geq 0}$ be a sequence  in the torus  converging to a point  $\tilde{\theta}$. The  kernel  of the sequence of domains  $\{\mathcal{U}_{\theta_{n}}\}$ is  $\mathcal{U}_{\tilde{\theta}}$. In particular the sequence  $\{\mathcal{U}_{\theta_n}\}$ kernel converges to $\mathcal{U}_{\tilde{\theta}}$.
 \end{lema}
 $\Pru$ Let  $\nabla$ be the kernel of the sequence $\{\mathcal{U}_{\theta_n}\}$. We first show   $\mathcal{U}_{\tilde{\theta}}\subset \nabla$. Let  $z$ be in $\mathcal{U}_{\tilde{\theta}}$. Let  $n^*\geq 2$ such that  $z$ belongs to  $\Omega_{n^*}(\mathcal{U}_{\tilde{\theta}})$. The set $\mathcal{U}$ is open,  thus there exists an open interval  $I_{\tilde{\theta}}$ around  $\tilde{\theta}$ such that if  $s\in I_{\tilde{\theta}}$ then  $\overline{\Omega_{n^*}(\mathcal{U}_{\tilde{\theta}})}\subset \mathcal{U}_s$. One concludes that the domain  $\Omega_{n^*}(\mathcal{U}_{\tilde{\theta}})$ is contained in  $\mathcal{U}_{\theta_n}$ for $n$ large enough. This shows that  $z$ belongs  to the kernel $\nabla$ of the domains  $\mathcal{U}_{\theta_n}$, i.e.
 \begin{equation}\label{introUennabla}
 \mathcal{U}_{\tilde{\theta}}\subset \nabla .
 \end{equation} 
In particular we see that the kernel  $\nabla$ is not reduced to the complex origin. The kernel is bounded since the domains  $\mathcal{U}_{\theta}$ are  uniformly bounded. Therefore $\nabla$ is a domain. Its conformal radius  verifies
 \begin{equation}\label{introRnablamayorRF}
 R(\nabla)\geq R({\tilde{\theta}})=R^{\mathcal{U}}_F.
 \end{equation} 
 Now we show  that $\nabla=\mathcal{U}_{\tilde{\theta}}$. Let $\eps>0$ and let $\tilde{n}\geq 2$ such that 
 \begin{equation*}
 R(\nabla)>R(\Omega_{\tilde{n}}(\nabla))>R(\nabla)-\eps .
 \end{equation*}
 From the definition of the kernel, for all  $z$ which belongs to the compact set  $\overline{\Omega_{\tilde{n}}(\nabla)}$ there exists a domain  $O_z$ containing the point   $z$ and contained in the domains $\mathcal{U}_{\theta_n}$ for $n$ large enough. Considering  a finite cover of $\overline{\Omega_{\tilde{n}}(\nabla)}$ by domains  $O_z$, we see that for  $n$ large  enough the compact set  $\overline{\Omega_{\tilde{n}}(\nabla)}$ is contained in  $\mathcal{U}_{\theta_n}$ since each $O_z$ has this property. This implies that for $n$ large enough
 \begin{equation*}
 R^{\mathcal{U}}_F=R({\theta_n})>R(\Omega_{\tilde{n}}(\nabla))>R(\nabla)-\eps
 \end{equation*}
 and we conclude that the conformal radius  $R^{\mathcal{U}}_F\geq R(\nabla)$, since the last inequality  holds for all $\eps>0$. Using  (\ref{introRnablamayorRF}), (\ref{introUennabla}) we have  the equalities 
 \begin{equation*}
 R^{\mathcal{U}}_F=R(\nabla)\quad , \quad \mathcal{U}_{\tilde{\theta}}=\nabla .
 \end{equation*}
 Any sub-sequence from  $\{\theta_n\}$ also converges to  $\tilde{\theta}$, and therefore all sub-sequence of domains from  $\{\mathcal{U}_{\theta_n}\}$ has the same kernel  $\nabla=\mathcal{U}_{\tilde{\theta}}\quad_{\blacksquare}$
 \paragraph{Proof of proposition \ref{prop2}. } 
 Let  $\{h_{\theta}:\D\to \mathcal{U}_{\theta}\}_{\theta\in \Td}$ be the family of uniformizing functions associated to the open tube  $\mathcal{U}$, with $h'_{\theta}(0)>0$, and let  $\Phi:\T\times \D \to \mathcal{U}$ be the bijective function defined by
 \begin{equation}
 (\theta,z)\longmapsto \big(\theta,h_{\theta}(z)\big).
 \end{equation}
From  lemma \ref{stable.232} and theorem \ref{caratheodory}, $\Phi$ is an homeomorphism between  $\Td\times \D$ and $\mathcal{U}$.
  Let  $\theta$ in $\Td$. The function $l_{\theta}:\D\to \D$ defined by
 \begin{equation}
 l_{\theta}(z)=\Pi_2\Big(\Phi^{-1}\circ F\circ \Phi(\theta,z)\Big)=h_{\theta+\alpha}^{-1}\Big(f_{\theta}\big(h_{\theta}(z)\big)\Big)
 \end{equation}
 is an holomorphic function, is univalent and verifies
 \begin{eqnarray*}
 l_{\theta}(0)&=&0\\
 l_{\theta}'(0)&=&\big(h_{\theta+\alpha}'(0)\big)^{-1}\rho_1(\theta)h_{\theta}'(0)\\
 &=&(R^{\mathcal{U}}_F)^{-1}\rho_1(\theta)R^{\mathcal{U}}_F\\
 &=&\rho_1(\theta).
 \end{eqnarray*}
  Since  $|\rho_1(\theta)|=1$, the Schwartz Lemma implies $l_{\theta}(z)=\rho_1(\theta)z$, and the result follows $\quad_{\blacksquare}$
\subsection{Proof of theorem \ref{teo1}}
We follow the original proof of R.P\'erez-Marco. We will need a finite version of the  P\'erez-Marco's theorem which will allows us to show the theorem for a fibered holomorphic dynamics having a rational (periodic) rotation on the base; an argument of continuity for the Hausdorff distance on compact sets will conclude the proof. Let $p\in \N$. We say that a $p-tuple$ of holomorphic transformations $(f_1,f_2,\dots,f_p)$  is a \emph{good chain} for the $p-tuple$ of sets $(U_1,U_2,\dots,U_p)$ if for every $i\in \{1,\dots,p\}$ we have
 \begin{itemize}
 \item[$i)$ ] $U_i$ is an open Jordan domain (the interior of a Jordan curve) containing the origin.
 \item[$ii)$ ] $f_i$ is a local holomorphic diffeomorphism fixing the origin, $\partial_zf_i(0)=e^{2\pi i \beta_i}$, $\beta_i\in \T$.
 \item[$iii)$ ] $f_i$ is defined and univalent in a neighborhood of  $\overline{U_i}$ and $f_i^{-1}$ is defined and univalent in a neighborhood of  $\overline{U_{i+1}}$, with $U_{p+1}=U_1$.  
 \end{itemize} 
 In order to fix ideas we will  perform the proofs for the case  $p=2$, even if the arguments also hold  in the general case.  Let  $(f_1,f_2)$ be a good chain for the pair  of sets $(U_1,U_2)$. We put  $f=f_2\circ f_1$. The function $f$ is a local  holomorphic diffeomorphism  with a fixed point in the origin and multiplicator equals to $e^{2\pi i (\beta_1+\beta_2)}$. 
 \begin{rema}
 Since the function  $f=f_2\circ f_1$ is an holomorphic function, it may happen that its maximal domain of definition is bigger than  (a neighborhood of) $\overline{U}_1\cap f^{-1}_1\big(\overline{U}_2\big)$, the set where the composition is \emph{a priori} defined, because  of analytic continuations. In the arguments that follow, we will consider the compositions as being defined only on (a neighborhood of) $\overline{U}_1\cap f^{-1}_1\big(\overline{U}_2\big)$.  
 \end{rema}
 
 It is a particular case of a simple but important result of Kerekjarto (see \cite{KERE23}, \cite{LEYO97}) that any connected component of ${U}_1\cap f^{-1}_1\big({U}_2\big)$ is a Jordan domain. We will denote by $U_{12}$ (resp. $U_{21}$) the connected component of ${U}_1\cap f^{-1}_1\big({U}_2\big)$ (resp. ${U}_2\cap f^{-1}_2\big({U}_2\big)$) which contains $0$.
    \paragraph{Fundamental construction. } Let  $K_1\subset \overline{U_1}$ and $K_2\subset\overline{U_2}$ be compact, connected  and full sets containing the origin and verifying
    \begin{equation*}
    f_1(K_1)=K_2=f^{-1}_2(K_1)\ , \ f_2(K_2)=K_1=f^{-1}_1(K_2)  
    \end{equation*}
    We will associate to these sets an analytical diffeomorphism of the circle. Indeed, let  $h_1:\overline{\C}\setminus \overline{\D}\to \overline{\C}\setminus K_1$ be an uniformization with $h_1(\infty)=\infty$ and $h_2$ be defined in an anologous way for $K_2$. The composition 
    \begin{equation*}
    h_1^{-1}\circ f_2\circ h_2\circ h_2^{-1}\circ f_1 \circ h_1=h_1^{-1}\circ f\circ h_1
    \end{equation*}  
    is an holomorphic diffeomorphism  $g$ from  an open annular set contained in  $\C\setminus \overline{\D}$, which has  $\T=\partial \D$ as one of its boundaries, onto a similar domain. The Carath\'eodory's extension theorem tell us that  this diffeomorphism extends continuously to a homeomorphism up to the boundary $\T$. The Schwartz's reflection principle allows us to extend   $g$ to an annular neighborhood  of   $\T$, producing in this manner an holomorphic diffeomorphism defined on this set. This diffeomorphism preserves the circle and then results  in an analytic circle diffeomorphism. 
%      \begin{figure}\label{construccionfundamental}
% \centering
%\includegraphics[width=100mm]{contruccionfundamental.png}
%\caption{Esquise de la construction fondamentale}
%\end{figure}

    \paragraph{The linearizable case. } Assume  $\beta=\beta_1+\beta_2$  is a real number verifying a diophantine condition. The Siegel linarization theorem  (see \cite{SIEG42}) implies $f$ is linearizable in an open neighborhood of the fixed point $z=0$. Let  $S(f)\subset U_{12}$ be the maximal linearization domain for  $f$  in $U_{12}$, let  $K_1=\widehat{\overline{S(f)}}$ be the compact set obtained by filling the closure of  $S(f)$ and let  $K_2=f_1(K_1)$. As $f_1$ is an homeomorphism on a neighborhood of $K_1$, $K_2$ is the filling of the closure of $f_1\big(S(f)\big)$. We apply the fundamental construction to  $K_1$ and $K_2$ in order to get an analytic circle diffeomorphism $g$.
    \begin{lema}
    Under the above hypothesis the analytic circle  diffeomorphism $g$ has rotation number   $\varrho(g)=\beta$.
    \end{lema}  
    $\Pru $ We follow an argument from \cite{PERE97}. Let be  $\{\Omega^1_n\}_{n\geq 0}$ a sequence of regions in $S(f)$ such that
    \begin{itemize}
    \item[\emph{i)} ] $\overline{\Omega^1_n}\subset \Omega_{n+1}^1$ for every $n\geq 0$.
    \item[\emph{ii)} ] $\bigcup_{n\geq 0}\Omega_n^1=S(f$).
    \item[\emph{iii)} ] $\overline{\Omega^1_n}$ is completely invariant under  $f$.
    \end{itemize} 
    Such a sequence may be obtained as the image of radius increasing radii discs by the linearizing application of  $f$. We also define a sequence  $\{\Omega_n^2\}_{n\geq 0}$ by $\Omega_n^2=f_1(\Omega_n^1)$ which verify analogous exhausting properties for the set  $f_1\big(S(f)\big)$. Let $\{g_n\}_{n\geq 0}$ be the sequence of circle analytic diffeomorphisms resulting from the application of the fundamental construction to  $\overline{\Omega_n^1}$ and $\overline{\Omega_n^2}$. Since $g_n$ is linearizable for all $n\geq 0$ (because of $\overline{\Omega_n^1}\subset S(f)$) we have $\varrho(g_n)=\beta$ for every $n\geq 0$. Let $h_i^n$ be the uniformizations  $\overline{\C}\setminus \overline{\Omega_n^i}$ for 
    all $n\geq 0$ and $i$ in $\{1,2\}$. From the definition of   $K_1$, the sequence $\{\overline{\C}\setminus\overline{\Omega_n^i}\}_{n\geq 0}$ is kernel convergent to $ \nabla \big(\{\overline{\C}\setminus\overline{\Omega_n^i}\}_{n\geq 0}\big)=\overline{\C}\setminus K_i$, 
    for $i$ in $\{1,2\}$. This implies that the holomorphic diffeomorphisms $g_n$ are defined simultaneously in some annular neighborhood  $A$ of $\T$ for every  $n$ large enough. Let $C\subset A$ be  a curve, which is  homotopic and exterior to $\T$. Let $C'$ be its reflection which respect to  $\T$. The Kernel Carathodory's theorem gives
    \begin{equation*}
    g_n \big |_{C}=\big (h_1^n\big )^{-1} \circ f \circ h_1^n \big |_{C}\stackrel{n\to \infty}{\longrightarrow} g\big|_{C}=h_1 \circ f \circ h_1
        \end{equation*}
    uniformly on $C$. We have the same result on  $C'$, and by the maximum principle applied to the compact annulus     ${A'}$,  bounded by  $C$ and $C'$, we have  $g_n\big|_{A'}\to g\big|_{A'}$ uniformly. Therefore we  get that $\varrho(g)=\beta$ since $\varrho(g_n)=\beta$ for all $n\geq 0\quad_{\blacksquare}$
    \\
    
    The global linearization theorem for analytic circle diffeomorphisms  (see \cite{HERM79}, \cite{YOCC99}) implies thus that $g$ is linearizable in an annular neighborhood of  $\T$. 
    If $S(f)$ was relatively compact in $U_{12}$, we would get  a larger linearization domain in $U_{12}$ by gluing $K_1$ to the image of $h_1$ of a small $g-$invariant neighborhood of $\T$. As $S(f)$ was the maximal linearization domain, we must have $\partial S(f)$ meets $\partial U_{12}\subset \big(\partial U_1\cup f^{-1}_1(\partial U_2)$. Thus, either $K_1$ meets $\partial U_1$ or $K_2$ meets $\partial U_2$ (or both). We have showed     
    \begin{prop}\label{introerizofinito}
    Let $(f_1, f_2, \dots ,f_p)$ be a good chain for the open sets  $(U_1,\dots,U_p )$. Assume the sum of the rotation numbers $\beta=\sum_i \beta_i$  verifies a diophantine condition and let  $S(f)\subset U_1$ be the maximal linearisation domain for $f=f_p\circ\dots\circ f_1$ in $U_1\cap f_1^{-1}(U_2)\cap \dots \cap \big(f_{p-1}\circ\dots\circ f_1\big)^{-1}(U_p)$. We consider the compact, connected full sets containing the origin given by $K_1=\widehat{\overline{S(f)}}, K_2=f_1(K_1), \dots, K_p=f_{p-1}(K_{p-1})$. There exists an index  $i$ in $\{1,2,\dots,p\}$ such that  $K_i\cap \partial U_i\neq \emptyset$. 
    \end{prop}
    In the following paragraphs we want to  eliminate the arithmetical hypothesis on $f$. We will show indeed
    \begin{prop}\label{introerizonolineal}
      Let $(f_1, f_2, \dots ,f_p)$ be a good chain for the open sets   $(U_1,\dots,U_p )$. There exists compact, connected full sets  containing the origin $K_1\subset \overline{U_1}, K_2\subset \overline{U_2},\dots, K_p\subset \overline{U_p}$ verifying
      \begin{equation*}
      f_i(K_i)=K_{i+1}\quad , \quad f_{i}^{-1}(K_{i+1})=K_i 
      \end{equation*} for all  $i$ in $\{1,2,\dots,p\}$, and an index $i^*$  such  that 
      $K_{i^*}\cap \partial U_{i^*}\neq \emptyset$. 
    \end{prop}  
%      \begin{figure}\label{erizofibradofinito}
% \centering
%\includegraphics[width=100mm]{erizofibrado.png}
%\caption{La Proposition \ref{introerizonolineal}}
%\end{figure}
    
    The proof of this proposition is obtained by showing that the conclusion is a dense and closed property in an adequate space. We will need some topology about compact sets.
    \paragraph{Hausdorff distance. }  Let $(X,d)$ be a compact metric space. Let $\eps>0$ and $A\subset X$. We denote by $V_{\eps}(A)$ the $\eps-$neighborhood of $A$ in $X$
    \begin{equation*}
    V_{\eps}(A)=\big\{ x\in X\ \big|\  d(x,A)<\eps\big\}
    \end{equation*}
    The Hausdorff distance between two non-empty compact sets $K_1,K_2$ contained in  $X$ is defined by 
    \begin{equation*}
    d_H(K_1,K_2)=\inf\big\{\eps>0\  \big| \  K_1\subset V_{\eps}(K_2)\quad and\quad K_2\subset V_{\eps}(K_1) \big\}
    \end{equation*}
    The  $d_H$ is a distance on the space $\mathcal{K}(X)$ of compact non-empty subsets of $X$. Moreover, $(\mathcal{K}(X),d_H)$ is a compact metric space. The space $\mathcal{K}_c(X)$ of connected,  non-empty compact sets is a compact subspace of  $\mathcal{K}(X)$. Let  $(Y,d')$ be another compact metric space. We denote by $C^0(X,Y)$ the space of continuous functions from $X$ to  $Y$ and we endow it with    the topology of the uniform convergence.
    \begin{lema}[see \cite{PERE97}]\label{Hcontinuo}
    The application
    \begin{eqnarray*}
    C^0(X,Y)\times \mathcal{K}(X)&\longrightarrow& \mathcal{K}(Y)\\
    \big(f,K\big)&\longmapsto& f(K)
    \end{eqnarray*}
    is continuous $\quad_{\blacksquare}$
    \end{lema}
    \begin{lema}\label{lemaA}
    Let $\{A^n\}_{n\in \N}$ be a sequence of compact sets of $X$ and $A$ be a compact set such that $A^n\to A$ in the Hausdorff distance. Let $x\in X$. There exists a sequence $x_n\in A^n$ such that $x_n\to x$ iff $x\in A \quad_{\blacksquare}$
    \end{lema}
   \begin{lema}\label{lemafull}
   Let $A\subset \T\times \C$ be a compact set, whose fibers are connected. Let $\widehat{A}_{\theta}$ be the filling of the fiber $A_{\theta}$, for every $\theta\in \T$. The set $\widehat{A}=\bigcup_{\theta\in \T}\{\theta\}\times \widehat{A}_{\theta}$ is compact. 
   \end{lema}
   $\Pru $ We must to show that $\widehat{A}$ is closed. Otherwise, suppose $(\theta_n,z_n)\in \{\theta_n\}\times \widehat{A}_{\theta_n}$ is a sequence converging to a point $(\theta,z)\notin \{\theta\}\times \widehat{A}_{\theta}$. The set $\widehat{A}_{\theta}$ is a compact full set, therefore there exists a continuous path $\gamma\subset \C\cup\{\infty\}$ which joints $z$ to $\infty$ and that does not touches $\widehat{A}_{\theta}$. Define a sequence of continuous paths $\gamma_n$ joining $z_n$ to $\infty$ by $\gamma_n=[z_n,z]\cup \gamma$. We have $\gamma_n\to \gamma$ in the Hausdorff topology on $\C\cup\{\infty\}$. Since $z_n\in \widehat{A}_{\theta_{n}}$, for every $n\in \N$ there exists a point $\tilde{z}_n\in A_{\theta_n}\cap\gamma_n$. Choose a sub-sequence $\{\tilde{z}_{n_i}\}_{i\geq 0}$ and a point $\tilde{z}$ such that $\tilde{z}_{n_i}\to \tilde{z}$. Thus $\tilde{z}\in A_{\theta}\cap \gamma$, but $\gamma\cap \widehat{A}_{\theta}=\emptyset$, a contradiction $\quad_{\blacksquare}$ 
   \\
     
    Let  $\mathcal{F}$ be the set  of all the  $p-$tuples $(f_1,\dots,f_p)$ verifying the hypothesis of proposition \ref{introerizonolineal} for the open sets $(U_1,\dots,U_p)$. We embed $\mathcal{F}$ in the product space $$C(\overline{U_1},\C)\times C(\overline{U_2},\C)\times C(\overline{U_2},\C)\times C(\overline{U_3},\C)\times \dots\times C(\overline{U_p},\C)\times C(\overline{U_1},\C)$$ by the inclusion
 \begin{equation*}
 (f_1,\dots,f_p)\longmapsto (f_1,f_1^{-1},f_2,f_2^{-1},\dots,f_p,f_p^{-1}).
 \end{equation*}
 The topology on $\mathcal{F}$ is the induced topology from the uniform convergence topology in the product. 
     \begin{lema}
    The set of $p-$tuples in $\mathcal{F}$ verifying the conclusions of proposition \ref{introerizonolineal} is closed.
    \end{lema}  
    $\Pru $ Let  $\{(f_1^n,\dots,f_p^n)\}_{n\geq 0 }$ be a sequence of $p-$tuples in $\mathcal{F}$ verifying  the conclusions of proposition \ref{introerizonolineal} and converging to the  $p-$tuple $(f_1,\dots,f_p)$ in $\mathcal{F}$. So, there exists a sequence of connected full compact sets  containing the origin $\{K_i^n\}_{n\geq 0}$ in $\mathcal{K}_c(\overline{U_i})$, $i\in \{1,\dots,p\}$, satisfying 
    \begin{eqnarray*}
    f_i^n(K_i^n)&=&K_{i+1}^n\\
    \big(f_i^n\big)^{-1}(K_{i+1}^n)&=&K_i^n.
    \end{eqnarray*}
    Moreover, for all $n\geq 0$ one has $K_i^n\cap\partial U_{i}\neq \emptyset$ for at least one integer $i\in \{1,\dots, p\}$. Since $\mathcal{K}_c(\overline{U_i})$ is compact for every $i$, passing to a subsequence if necessary, we get compact connected sets  $\tilde{K}_i$ in $\mathcal{K}_c(\overline{U_i})$ for $i\in \{1,\dots,p\}$ verifying for every $i$ that
    \begin{equation*}
    K_i^n\stackrel{n\to \infty}{\longrightarrow} \tilde{K}_i
    \end{equation*} 
    in the Hausdorff distance. Furthermore, there exists an index $i^*$ in $\{1,\dots,p\}$ such that $K_{i^*}^n\cap \partial{U_{i^*}}\neq \emptyset$ for infinitely many  values of $n$ and therefore $\tilde{K}_{i^*}\cap \partial{U_{i^*}}\neq \emptyset$. In the limit we have $f_i(\tilde{K}_i)=\tilde{K}_{i+1}$ and $f_i^{-1}(\tilde{K}_{i+1})=\tilde{K}_i$ for all $i$ in $\{1,\dots,p \}$, mod $p$. The set $\tilde{K}_1$ must be contained in the connected component of the origin in    $\overline{U}_1\cap f_1^{-1}(\overline{U}_2)\cap \dots \cap \big(f_{p-1}\circ\dots\circ f_1\big)^{-1}(\overline{U}_p)$, which is the closure of a Jordan domain, and similarly for $\tilde{K}_{i}$. Therefore, the filled-in sets $K_i=\widehat{\tilde{K}_i}$ verify the conclusions of proposition \ref{introerizonolineal} for the $p-$tuple $(f_1,\dots, f_p)\quad_{\blacksquare}$ 
    \\
    
    Now we can conclude the proof of proposition  \ref{introerizonolineal}.  If  $(f_1,\dots, f_p)$ belongs to $\mathcal{F}$, and has multipliers $(e^{2\pi i\beta_1},\dots, e^{2\pi i \beta_p})$, for  $i$ in $\{1,\dots,p\}$, there exists a sequence of real numbers  $\{\eta_n^i\}_{n\geq 0}$  converging to zero such that the  $\sum_{i}(\beta_i+\eta_n^i)$ is a real number verifying a diophantine condition for all  $n\geq 0$. The $p-$tuple 
    \begin{equation*}
    \big(e^{2\pi i \eta_n^1}f_1,\dots, e^{2\pi i \eta_n^p}f_p\big)
    \end{equation*}
    belongs to $\mathcal{F}$ for $n$ large enough and converges  to $(f_1,\dots,f_p)$. Moreover, they verify the conclusions of proposition  \ref{introerizonolineal}  by proposition \ref{introerizofinito}. The precedent lemma allows us to conclude $\quad_{\blacksquare}$ 
    \begin{rema}
    Proposition \ref{introerizonolineal} can be obtained in an alternative way as a direct application of the P\'erez-Marco's theorem (which corresponds to the case $p=1$). Indeed, suppose $p=2$. Under the hypothesis of the proposition, the composition $f=f_2\circ f_1$ and its inverse $f^{-1}=f_1^{-1}\circ f_2^{-1}$ are well defined and univalent in a neighborhood of the set
    \begin{equation*}
    \mathcal{V}=\bigodot\Big\{ U_1\cap f^{-1}_1(U_2)\cap \big(f_2^{-1}\big)^{-1}(U_2)\Big\},
    \end{equation*}
    where $\bigodot$ denote the connected component containing the origin. Then $\mathcal{V}$ is a Jordan domain  and the P\'erez-Marco's theorem applies. We have thus a connected, compact set $K$ containing the origin, completely invariant by $f$, touching the boundary of $\mathcal{V}$. Suppose for instance that $K\cap \partial f^{-1}_1(U_2)\neq \emptyset$. The sets $K_1=K$,   $K_2=f_1(K)$ satisfy the inclusions properties and $K_2\cap \partial U_2\neq \emptyset\quad_{\blacksquare}$
\end{rema}    

       We will show theorem \ref{teo1} as a consequence of proposition \ref{introerizonolineal}.  
%    \begin{theo}\label{introteoorbitaestable}
%    Let  $F$ be a \emph{fhd} having an invariant indifferent curve $u$. Let $\mathcal{U}$ be an open tube containing $u$, whose fibers are Jordan domains. We suppose $\overline{\mathcal{U}_{\theta}}$ depends continuously on $\theta$ and  $F$ and $F^{-1}$ define \emph{fhd}'s, which are  injective in some neighborhood of  $\overline{\mathcal{U}}$. Then there exist a compact connected set  $K\subset \overline{\mathcal{U}}$ such that
%    \begin{itemize}
%    \item[i) ] $graph(u)\subset K$.
%    \item[ii) ] $F(K)=F^{-1}(K)=K$, that is, $K$ is a completely invariant set for  $F$.
%    \item[iii) ] $K\cap \partial\mathcal{U}\neq\emptyset$.
%    \end{itemize} 
%    \end{theo}
%      \begin{figure}\label{erizototal}
% \centering
%\includegraphics[width=100mm]{erizototal.png}
%\caption{}
%\end{figure}
The actual proof is done in a very analogous way to the proof of  proposition \ref{introerizonolineal} (and also very analogous to the original proof of    P\'erez-Marco's theorem). We will allow  \emph{fhd}'s having a base rotation vector different of $\alpha$ on the torus $\Td$ (in fact, a rational rotation vector).   We recall that the invariant curve is supposed to be the zero section curve $u=\{z\equiv 0\}_{\Td}$. 
  
  Let $\mathcal{F}_{\mathcal{U}}$ be the set of all \emph{fhd} $F$ (with free rotation vector over $\Td$),  having the zero section as an invariant curve and  such that  $F$, $F^{-1}$ define injective \emph{fhd}'s in a neighborhood of $\overline{\mathcal{U}}$.  Thus, if a  \emph{fhd} $F$ belongs to  $\mathcal{F}_{\mathcal{U}}$ then it has the form
  \begin{equation}\label{intronotacionhabitualfin}
  F(\theta,z)=\big(\theta+\alpha_F,\rho_1(\theta)z+z^2\rho(\theta,z)\big),
  \end{equation}
  where $\alpha_F$ is a  vector in $\Td$, $\rho_1$ is a non-vanishing continuous function 
  and $\rho$ is a continuous function holomorphic in each fiber. 
  We endow  this set with the topology of the uniform convergence of $F$ and $F^{-1}$ on $\overline{\mathcal{U}}$. 
    \begin{lema}
  The set of \emph{fhd}'s in $\mathcal{F}_{\mathcal{U}}$ verifying the conclusions of theorem \ref{teo1} is closed.
  \end{lema}
  $\Pru $ Let $\{F_i\}_{i\geq 0}$ be a sequence in $\mathcal{F}_\mathcal{U}$ whose elements verify the conclusions of  theorem \ref{teo1}. We suppose the sequence converges to the \emph{fhd} $F$. Let $\{K^i\}_{i\geq 0}$ be an associated  sequence of compact connected sets  verifying the conclusions of the theorem. For every  $i\geq 0$ the invariant curve is included in the compact set  $K^i$, and there exist a point    $(\theta_i,z_i)$ in  $K^i\cap \partial \mathcal{U}$. The set $ \partial \mathcal{U}$ is a compact set and the space $\mathcal{K}_c( \overline{\mathcal{U}})$, formed by all connected compact sets of  $ \overline{\mathcal{U}}$, is also compact. Considering  a subsequence if necessary,  there exists a connected compact set  $K\subset\overline{\mathcal{U}}$ containing the invariant curve and a point  $(\tilde{\theta},\tilde{z})$ such that 
  \begin{eqnarray}\label{HdistanciaK}
  d_H(K^i, K)&<&\frac{1}{i}\\
  (\theta_i,z_i)&\stackrel{i\to \infty}{\longrightarrow}& (\tilde{\theta},\tilde{z}) \in K\cap \partial\mathcal{U}.
    \end{eqnarray}
 Lemma  \ref{Hcontinuo} implies that $K$ is completely invariant. Let's show the fibers $K_{\theta}$ are connected compact sets. We will use the following notation: let  $\tilde{\theta}\in \Td$, $\eps>0$ and $A\subset \overline{\mathcal{U}}$ a compact set. We define a compact set $A_{[\tilde{\theta};\eps]}$ by 
 \begin{equation*}
 A_{[\tilde{\theta};\eps]}=\bigcup_{|\theta-\tilde{\theta}|\leq \eps}\{\theta\}\times A_{\theta}.
 \end{equation*}
  Let $\theta\in \Td$. Inequality (\ref{HdistanciaK}) implies
  \begin{eqnarray}\label{iestrella}
  \{\theta\}\times K_{\theta}&\subset& V_{\frac{1}{i}}\big(K^i_{[\theta;\frac{1}{i}]}\big)\\ \label{iiestrella}
  K^i_{[\theta;\frac{1}{i}]}&\subset& V_{\frac{1}{i}}\big(K_{[\theta;\frac{2}{i}]}\big).
  \end{eqnarray}
  The sets $K^i_{[\theta;\frac{1}{i}]}$ are compact connected sets, therefore there exists a subsequence $\{K^{i_n}_{[\theta;\frac{1}{i_n}]}\}_{n\geq 0}$ converging in the Hausdorff topology to a connected compact set $\{\theta\}\times \tilde{K}_{\theta}$. We will show that $K_{\theta}=\tilde{K}_{\theta}$.  From (\ref{iestrella}), we have that $K_{\theta}\subset \tilde{K}_{\theta}$. Conversely, let $z\in \tilde{K}_{\theta}$. There exists a sequence $(\theta_{i_n},z_{i_n})\in K_{[\theta_{i_n};\frac{1}{i_n}]}^{i_n}$ such that $\theta_{i_n}\to \theta$ and $z_{i_n}\to z$. From (\ref{iiestrella}) we obtain a sequence $(\tilde{\theta}_{n},\tilde{z}_n)\in K_{[\theta;\frac{2}{i_n}]}$ such that $\big|(\tilde{\theta}_{n},\tilde{z}_n)-(\theta_{i_n},z_{i_n})\big|<\frac{1}{i_n}$. Therefore $(\tilde{\theta}_n,\tilde{z}_n)$ converges to $(\theta,z)$ and $z\in K_{\theta}$. Finally,  the compact set  $\widehat{K}$, obtained  filling  the fibers of $K$ as defined in lemma \ref{lemafull}, is still completely invariant and verifies all the conclusions of theorem \ref{teo1}  $\quad_{\blacksquare}$
  \\
  
  The proof of theorem  \ref{teo1} follows from
  \begin{lema}
  Let $F$ be a  \emph{fhd} verifying the hypothesis of theorem \ref{teo1}. There exists a sequence in $\mathcal{F}_{\mathcal{U}}$ converging to $F$, such that its elements verify  the conclusions of theorem  \ref{teo1}.
  \end{lema}
  $\Pru $ We write  $F$ as in  (\ref{intronotacionhabitualfin}). Let $\{l_n\}_{n\geq 0}$ be a sequence of zero-mean trigonometric polynomials, converging uniformly to the function  $\log |\rho_1|$ (F\'ejer's theorem).   Since $\alpha$ is an rationally independent vector we can pick a sequence of vectors $\alpha_n$ with rational components  $(p_n/q_n)_{1,\dots,d}$ converging to  $\alpha$. Furthermore, we may assume that this sequence verifies that   the denominators  $(q_n)_j$, $j\in \{1,\dots,d\}$ are different prime numbers bigger than the degree of the trigonometric polynomial $l_n$, for every $n\geq 0$. We put 
%   On \'ecrit la s\'erie de Fourier de $\log \big|\rho_1\big|$ par
%  \begin{equation*}
%  \log \big|\rho_1(\theta)\big|=\sum_{j\in \Z\setminus\{0\}}\hat{\rho}_1(j)e^{2\pi i j\theta} 
%  \end{equation*}
% Pour $n$ dans $\N$ nous possons 
%  \begin{equation}
%  \log \big|\rho_1\big|\Big|_n(\theta)=\sum_{0<|j|\leq n}\hat{\rho}_1(j)e^{2\pi ij \theta}
%  \end{equation}
%  la troncature au niveau $n$ du $\log \big|\rho_1\big|$, et finalement la fonction 
  \begin{equation}
  \rho_1^n=\frac{\rho_1}{\big|\rho_1\big|}e^{l_n}.
  \end{equation}
   %Comme $\rho_1$ ne s'anulle pas, le module $\big|\rho_1\big|$ est de classe $C^1$ et la convergence de la s\'erie de Fourier est uniforme dans $\T$, donc $\rho_1^n$
    These functions converge uniformly to $\rho_1$. We define a sequence of  \emph{fhd}'s $\{F_n\}_{n> 0}$ by 
   \begin{equation}
   F_n(\theta,z)=\Big(\theta+\alpha_n,\rho_1^n(\theta)z+z^2\rho(\theta,z)\Big).
   \end{equation}
   For every  $n$ large enough the transformation  $F_n$ belongs  to $\mathcal{F}_{\mathcal{U}}$ and moreover $F_n$ converges to $F$.We can solve the cohomological equation 
     \begin{equation}
    \tilde{u_1}^n(\theta+\alpha_n)-\tilde{u_1}^n(\theta)=l_n(\theta)
    \end{equation}
    since the conditions over the denominators    avoids to fall in a resonance using the Fourier coefficients method.  We put $u_1^n=\exp(\tilde{u_1}^n)$. This function satisfies 
    \begin{equation}
    \frac{u_1^n(\theta+\alpha_n)}{u_1^n(\theta)}=e^{l_n(\theta)}.
    \end{equation} 
    The change of coordinates   $\tilde{F}_n=H_n^{-1}\circ F_n\circ H_n$, with 
    \begin{equation*}
    H_n(\theta,z)=\big(\theta,u_1^n(\theta)z\big)
    \end{equation*}
    gives us the following  normal form for $F_n$
    \begin{equation*}
    \tilde{F}_n(\theta,z)=\Big(\theta+\alpha_n,\tilde{\rho}_1^n(\theta)z+z^2\tilde{\rho}^n(\theta,z)\Big),
    \end{equation*}
    with $|\tilde{\rho}_1^n(\theta)|=1$ for every $\theta$ in $\T$. %%%%%% 
    We put $\mathcal{U}^n=H_n^{-1}\big(\mathcal{U}\big)$. Let $\tilde{\theta}$ in $\Td$ and $Q_n=\prod_{1}^d (q_n)_j$. We consider the holomorphic transformations   
    \begin{equation*}
    {\gamma}_i^n(z)=\tilde{f}_n\Big(\tilde{\theta}+i\alpha_n,z\Big)
    \end{equation*}
    for  $i$ in $\{0,\dots,Q_n-1\}$. Thus, the $Q_n-$tuple $(\gamma_0^n,\dots,\gamma_{Q_n-1}^n)$ is a good chain for the fibers  $\mathcal{U}^n_{\tilde{\theta}+i\alpha_n}$ respectively. Applying proposition \ref{introerizonolineal} we get connected, full and compact sets  $K_i^{n}\subset \mathcal{U}^n_{\tilde{\theta}+i\alpha_n}$, $i\in \{0,\dots, Q_n-1\}$, containing the origin and verifying
    \begin{eqnarray*}
   \gamma_i^n(K_i^n)&=&K_{i+1}^n\\
    (\gamma_{i+1}^n)^{-1}&=&K_i^n.
    \end{eqnarray*}
    We also have an index  $i=i_n$ in $\{0,\dots,Q_n-1\}$ such that  $K_i^n\cap \partial\mathcal{U}^n_{\tilde{\theta}+i\alpha_n}\neq \emptyset$. The set
    \begin{equation*}
    K^n=H_n\bigg(\{z\equiv 0\}_{\T}\bigcup \Big(\bigcup_{j=0}^{Q_n-1}\Big\{\tilde{\theta}+j\alpha_n\Big\}\times K_j^n\Big)\bigg)
    \end{equation*}
    verifies the conclusions of theorem \ref{teo1} for $F_n$ $\quad_{\blacksquare}$
      \\
      
     \begin{coro}
     Let $F$ be a  \emph{fhd} and $u$ be an indifferent invariant curve. Then there exists orbits, other to those situated on the curve, that are  in past and future completely contained in a neighborhood of the curve. In fact, actually there exists such an orbit in any tubular neighborhood of the invariant curve.
          \end{coro}
          $\Pru $ As usual we suppose  $u=\{z\equiv 0\}_{\Td}$ and $F(\theta,z)=\big(\tal,\rho_1(\theta)z+z^2\rho(\theta,z)\big)$, with $\rho_1(\theta)\neq 0$ for every  $\theta$ in $\Td$. There exists a constant  $C>1$ such that  $C>|\rho_1(\theta)|>C^{-1}$ for every $\theta$ in $\Td$. Applying  lemma  \ref{lemainversion} we get an uniform positive radius  $\bar{r}$ such that  $f_{\theta}$ and $f_{\theta}^{-1}$ are well defined and univalent in some neighborhood of   $\overline{D_r}$ for every $r\in (0,\bar{r}]$. Taking  $\mathcal{U}=\Td\times D_r$,  theorem \ref{teo1} applies $\quad_{\blacksquare}$ 
          \\
          
          $\mathbf{Question. }$ We had showed that the invariant continua $K$ has at least (and so infinitely many) one fiber that is different of the origin alone. One can ask if this property holds for every fiber. We conjecture that this is not always  the case and that there should exists examples where singular fibers (reduced to the origin alone) arise. 
\bibliographystyle{plain}
\bibliography{topotesisd.bib}
\end{document}